\documentclass[12pt]{amsart}

\usepackage{amssymb,amscd,amsthm,mathrsfs}
\usepackage[ansinew]{inputenc}
\usepackage[centertags]{amsmath}
\usepackage[all]{xypic}
\usepackage{graphicx}
\usepackage{color}

 \def\RR{{\mathbb R}} \def\SS{{\mathbb S}} \def\TT{{\mathbb T}}
   
 \def\ZZ{{\mathbb Z}}

    \def\cS{\mathcal{S}}

    \def\cW{\mathcal{W}}

\newtheorem*{teo*}{Theorem}

\newtheorem{teo}{Theorem}[section]

\newtheorem{quest}{Question}

\newtheorem{lema}[teo]{Lemma}
\newtheorem{prop}[teo]{Proposition}

\newcommand{\bi}{\begin{itemize}}
\newcommand{\ei}{\end{itemize}}

\theoremstyle{definition}

\theoremstyle{remark}
\newtheorem{obs}[teo]{Remark}

\newcommand{\dem}{\vspace{.05in}{\sc\noindent Proof.\,\,}}
\newcommand{\lqqd}{\par\hfill {$\Box$} \vspace*{.05in}}

\newcommand{\eps}{\varepsilon}

\newcommand{\en}{\subset}

\DeclareMathOperator{\Diff}{Diff}

\newcommand{\comment}[1]{}

\author[C. Bonatti]{Christian Bonatti}
\address{CNRS - IMB. UMR 5584.
\rm Universit\'e de Bourgogne, 21004 Dijon, France} \email{bonatti
"at" u-bourgogne.fr}

\author[R. Potrie]{Rafael Potrie}
\address{CMAT, Facultad de Ciencias, Universidad de la Rep\'ublica, Uruguay}
\urladdr{www.cmat.edu.uy/$\sim$rpotrie} \email{rpotrie "at"
cmat.edu.uy}

\title[Intermingled basins]{Many intermingled basins in dimension 3}
\thanks{R.P. was partially supported by CSIC group 618. }

\begin{document}

\begin{abstract}
We construct a diffeomorphism of $\TT^3$ admitting any finite or countable number of physical measures with intermingled basins. The examples are partially hyperbolic with splitting $T\TT^3 = E^{cs} \oplus E^u$ and can be made volume hyperbolic and topologically mixing. 

\bigskip


\medskip

\end{abstract}

\maketitle

\section{Introduction}\label{SectionIntroduccion}

In \cite{Kan} an example of disjoint attractors with intermingled basins of attraction was presented. The example posses two invariant measures,
which are generic for positive Lebesgue measure subsets of the manifold while these sets are \emph{intermingled} in the sense that their Lebesgue 
density points are dense in one another (and in the whole manifold). 
This seems to be a rare phenomena in the space of all diffeomorphisms (see for example \cite{UV,VY}) 
yet it possesses some sort of robustness (Kan's example is robust on manifolds with boundary). 

The example of \cite{Kan} admits only two such measures with intermingled basins, and it was recently proved in \cite{UV} that this is essentially 
the only way to present measures with intermingled basins for \emph{strongly partially hyperbolic diffeomorphisms} of 3-manifolds (see below for precise definitions). 

One may wonder if it is possible to construct examples where more than $2$ disjoint attractors have intermingled basins, 
and it turns out to be possible. In \cite{MW,VY} examples in dimensions $\geq 4$ were presented. 
In this paper we give examples in dimension 3, which are \emph{partially hyperbolic} (but not strongly partially hyperbolic), 
and have an arbitrarily large number of intermingled measures. 
This answers a problem posed in \cite[Remark 1.4]{MW}. 

\section{Precise statement of results}\label{section-statements}
\subsection{Definitions} Before we state the main result precisely, we  introduce some definitions.

\subsubsection{Partial hyperbolicity} Let $f: M \to M$ be a $C^1$-diffeomorphism and $\Lambda \en M$ a compact $f$-invariant subset.

Given a $Df$-invariant bundle $E \en T_\Lambda M$ we say it is \emph{uniformly contracted} if there exists $N>0$ such that for every $x\in \Lambda$ one has that $\|Df^N|_{E(x)} \| < \frac 1 2$. 
A bundle is \emph{uniformly expanded} if it is uniformly
contracted for $f^{-1}$. Given two $Df$-invariant subbundles $E, F
\en T_\Lambda M$ we say that $F$ \emph{dominates} $E$ if there
exists $N>0$ such that $\|Df^N|_{E(x)}\| \|Df^{-N}|_{F(f^N(x))}\|
< \frac 1 2$ for every $x\in \Lambda$.

We say that $\Lambda$ is:

\bi \item[--] a \emph{partially hyperbolic} set of $f$ if
$T_\Lambda M = E \oplus F$ is a $Df$-invariant splitting such that
$F$ dominates $E$ and either $E$ is uniformly contracted or $F$ is
uniformly expanded. 
\item[--] a \emph{volume hyperbolic} set of $f$ if it admits a $Df$-invariant splitting $T_\Lambda M = E \oplus F$ such that $F$ dominates $E$, the jacobian\footnote{Given a Riemannian metric one can measure the $k$-dimensional volume along a $k$-dimensional subbundle of $TM$; the \emph{jacobian} along a subspace measures how this volume changes.} of $Df$ along $E$ is uniformly contracted and the jacobian of $Df$ along $F$ is uniformly expanded. 
\item[--] a \emph{strongly partially
hyperbolic} set of $f$ if $T_\Lambda M = E^s \oplus E^c \oplus
E^u$ is a $Df$-invariant splitting where the bundles $E^s$ and
$E^u$ are non trivial, $E^s$ is uniformly contracted, $E^u$ is
uniformly expanded, $E^u$ dominates $E^s \oplus E^c$ and $E^c
\oplus E^u$ dominates $E^s$. \ei

See \cite[Appendix B]{BDV} for more detailed introduction and
relevant properties. We will use the fact that in a partially
hyperbolic set the strong bundles (the ones which are either
uniformly contracted or expanded) integrate uniquely into an
$f$-invariant lamination which will be called \emph{stable
lamination} or \emph{unstable lamination} depending on whether it
integrates $E^s$ or $E^u$.

\subsubsection{Ergodic attractors and intermingled basins}

Given an $f$-invariant ergodic measure $\mu$ we define its
\emph{basin of attraction} $B(\mu)$ as the set of points $x \in M$ such
that for every $\varphi \in C^0(M)$ one has that $$\frac 1 n
\sum_{i=0}^{n-1} \varphi (f^i(x)) \to \int \varphi d\mu.$$

We say that $\mu$ is an \emph{SRB measure} if the Lebesgue measure
of $B(\mu)$ is positive. If $\mu_1, \ldots, \mu_k$ are SRB
measures\footnote{Sometimes, these measures are called \emph{physical measures}, but in this paper it will not make a difference.}, we say that their basins are \emph{intermingled} if for
every $U$ open set such that $B(\mu_i) \cap U$ has positive
Lebesgue measure for some $i$ then $B(\mu_j) \cap U$ has positive
Lebesgue measure for every $j$.

If $\Lambda$ is a partially hyperbolic set saturated by unstable
manifolds we say that a $f$-invariant measure is a
\emph{Gibbs-u-state} if the disintegration along unstable
manifolds is absolutely continuous (see \cite[Chapter 11]{BDV}).
In some cases, one can ensure that $\Lambda$ has a unique
Gibbs-u-state which is also an SRB measure, one paradigmatic such
example is when $f$ is $C^2$ and the set is \emph{mostly
contracting} (see Section \ref{s.prelim} or \cite{BV,VY}).

We say that a compact $f$-invariant set $\Lambda$ is a
\emph{minimal attractor} (in the sense of Milnor \cite{Milnor}) if
the Lebesgue measure of the set of points $x$ such that $\omega(x)
\en \Lambda$ is positive and for every compact $f$-invariant
proper subset $\Lambda' \en \Lambda$ the set of points $x$ such
that $\omega(x) \en \Lambda'$ has zero Lebesgue measure. When an
SRB measure $\mu$ is also a Gibbs-u-state and it is \emph{mostly
contracting}, it is possible to show that the support of $\mu$ is
a minimal Milnor attractor.

\subsection{Statements}

The main result of this note can be stated as follows:

\begin{teo}\label{teo-main} For every $k\geq 2$, there exists a partially hyperbolic $C^\infty$ diffeomorphism $f: \TT^3 \to \TT^3$ with splitting $T\TT^3 = E^{cs} \oplus E^u$ admitting exactly $k$ SRB measures $\mu_1,\ldots, \mu_k$ such that:
\begin{itemize}
\item[--] The union of the basins of $\mu_i$ has full Lebesgue
measure in $\TT^3$. \item[--] All their basins are intermingled.
\item[--] Their supports are disjoint and if $\Lambda_i$ denotes
the support of $\mu_i$ then $\Lambda_i$ is strongly partially
hyperbolic and mostly contracting (in particular, the measures $\mu_i$ are Gibbs-u-states and their support is a minimal Milnor attractor).
\end{itemize}
\end{teo}

\begin{obs} In fact, it is possible to construct examples with \emph{countably} infinitely many intermingled basins by adapting the construction. We indicate the main steps in this generalisation in section \ref{sec-countable}. 
\end{obs}

\begin{obs}
Kan's example has only two intermingled basins and is transitive on $\TT^2 \times [0,1]$. It can be made transitive in $\TT^3$ by gluing two copies of the example and interchanging the dynamics; however, this forbids the example to be topologically mixing ($f^2$ ceases to be transitive). Our example can be done topologically mixing by using the liberties we have in the construction, in particular, we use the fact that the construction can be made \emph{volume hyperbolic} (this is done in section \ref{sec-volumehyperbolic} ). The construction for having a topologically mixing example is explained in section \ref{sec-recurrence}.
\end{obs}

\begin{obs} 
In dimension 2 there are results in the direction of the non-existence of disjoint attractors with intermingled basins. In particular, in \cite{HHTU} it is shown that a transitive surface $C^{1+\alpha}$-diffeomorphism cannot have more than one hyperbolic SRB-measure.  
\end{obs} 

As already said above, \cite{UV} shows that, among strongly partially hyperbolic diffeomorphisms of $\TT^3$, no more than two basins can be intermingled.
Theorem \ref{teo-main} was motivated by a question that arose in the conference ``International Conference on Dynamical Systems" held at IMPA in November 
of 2013 during a talk where R. Ures  was presenting the results of \cite{UV}. We thank S. Crovisier, L.Diaz, R. Ures, M. Viana and J. Yang for useful comments and remarks. 

\section{Results used and outline of the construction}\label{s.prelim}

The main result we will use has to do with the study of the basins of $u$-\emph{Gibbs measures} with negative Lyapunov exponents along the center-stable direction. The following result from \cite{VY} generalises the criteria in \cite{Kan} and depends on the analysis made in \cite{BV} (see also \cite[Chapter 11]{BDV}). Let us first introduce some notions. 

A partially hyperbolic diffeomorphism $f: M \to M$ with splitting $TM=E^{cs} \oplus E^u$ is said to be \emph{mostly contracting} if for every unstable disk $D$ (i.e. tangent to $E^u$) there exists a subset $D^0 \subset D$ with  positive Lebesgue measure (for the induced metric on the disk) such that for every $x\in D^0$ one has that:

\begin{equation}\label{eq:mostly} 
\limsup_{n\to \infty} \frac{1}{n} \log \|Df^n|_{E^{cs}(x)}\| < 0 \ .
\end{equation}

This property is equivalent to having that all $u$-Gibbs measures have negative Lyapunov exponents along $E^{cs}$. Indeed, for a $C^2$-partially hyperbolic diffeomorphism, and $D$  an unstable disk, almost every point in $D$ converges to some $u$-Gibbs state (see \cite[Section 11.2.1]{BDV}). If the Lyapunov exponents along $E^{cs}$ are all negative, then it follows that condition \eqref{eq:mostly} is verified. 

For a mostly contracting partially hyperbolic diffeomorphism, one defines a \emph{skeleton} (c.f. \cite{VY}) to be a finite set of hyperbolic periodic points $q_1, \ldots, q_k$ such that their stable dimension equals $\dim E^{cs}$ and such that the following conditions hold:

\begin{itemize}
\item[(S1)] for every $x\in M$ the unstable manifold of $x$ intersects the stable manifold of the orbit of $q_i$ for some $i$,
\item[(S2)] if $i\neq j$ then the stable manifold of the orbit of $q_i$ and the unstable manifold of the orbit of $q_j$ are disjoint. 
\end{itemize}

The following is a direct consequence of \cite[Theorem A]{VY}.

\begin{teo}\label{teo-basins}
Let $f: M \to M$ be a $C^2$ partially hyperbolic diffeomorphism with splitting $TM=E^{cs} \oplus E^u$ which is mostly contracting. Assume that $\cS=\{q_1, \ldots, q_k \}$ is a skeleton for $f$ and such that, for every $i$, the stable manifolds of the orbit of $q_i$ is dense, then: 
\begin{itemize}
\item $f$ has exactly $k$ SRB measures $\mu_1, \ldots, \mu_k$ with disjoint supports and such that the unstable manifold of the orbit of $q_i$ is dense in the support of $\mu_i$. 
\item the union of the basins of the $\mu_i$ cover Lebesgue almost every point in $M$ and their basins are all intermingled. 
\end{itemize}
\end{teo}

\begin{obs}\label{rem.mostlycontracting}
Under assumption (S1), every $u$-Gibbs state contains exactly one of the $q_i$. Therefore, one can show that $f$ is mostly contracting by showing that each $q_i$ is contained in the support of a $u$-Gibbs state whose Lyapunov exponent is negative. See \cite{BV,VY}. 
\end{obs}

Using this result let us give a brief outline on how to construct the examples announced by Theorem \ref{teo-main}. 

\medskip

We consider a fibered map $f_0$ over a linear Anosov diffeomorphism $A: \TT^2 \to \TT^2$ of the form $(x,t)\mapsto (Ax, g_x(t))$ where we choose $A$ conveniently so that it has many periodic points and $g_x : S^1 \to S^1$ are circle diffeomorphisms varying smoothly with respect to $x \in \TT^2$. 

We choose  also a finite number of points in $S^1$ that we declare to be fixed by $g_x$ for every $x\in \TT^2$. This way, we force the existence of several minimal $\cW^u$-saturated subsets of $\TT^3$ for $f_0$ which is partially hyperbolic with splitting $T\TT^3 =E^s \oplus E^c \oplus E^u$. We also require the $g_x$ to be such that:

\begin{itemize}
\item in every one of these minimal sets there will be at least one periodic point for which the center direction is contracting and its stable manifold accumulates in the consecutive minimal sets, 
\item the diffeomorphism $f_0$ is mostly contracting.
\end{itemize}

This can be thought of as a finite number of Kan's examples glued together. We also require the existence of a circle on which the dynamics moves in one direction (like saddle nodes in each fixed point, see figure \ref{figure1}). This is done in Section \ref{s.Kan}.

\smallskip 

Then we perform a modification of $f_0$ that can be though of as done in two steps: 

\begin{itemize}
\item First one creates $f_1$, a DA-like deformation along each invariant tori, in the fixed points associated with the saddle-node circle. (This is done in Section 5.)
\item Then, one pushes the dynamics along the "hole" created by the DA-deformation so that the points traverse from one region to the one which is above. (This is done in Section \ref{s.second}.)
\end{itemize}

One can then check (see Section \ref{SectionProof}) that the conditions of Theorem \ref{teo-basins} are satisfied and this is enough to complete the the proof of Theorem \ref{teo-main}.

Finally, in Section \ref{sec-countable} we explain how to extend the construction to obtain countably 
infinitely many intermingled basins, in Section \ref{sec-volumehyperbolic}  we show how to obtain a volume hyperbolic example and in Section \ref{sec-recurrence} we study some recurrence properties of the examples.



\section{Many Kan's examples}\label{s.Kan}

 Let us fix $k\geq 2$ an integer, and for simplicity of our presentation, we will assume that $k$ is divisible by $6$ 
 (see Remark~\ref{remark-odd} which expalins how to adapt the argument in
 the other cases). 
 We first construct a smooth diffeomorphism $f_0 \in \Diff^\infty (\TT^3)$ which concatenates $k$ ``Kan's examples'' together.

Let $A \in SL(2,\ZZ)$ be a matrix which defines a linear Anosov diffeomorphism $g_A: \TT^2 \to \TT^2$. For simplicity, we shall
assume that $g_A$ has at least  five fixed points. One can assume moreover that one has that vectors in $E^u_A$ are expanded by a factor of
$\lambda>5$ and in $E^s_A$ are contracted by a factor of $\lambda^{-1}< \frac 1 5$. We denote as $W^s_A$ and $W^u_A$ the
$g_A$-invariant stable and unstable foliations which are projections to $\TT^2$ of the bundles $E^s_A$ and $E^u_A$
respectively.

Consider $\TT^3 = \TT^2 \times S^1$ with coordinates $(x,t)$ where $x\in \TT^2$ and $t\in S^1$. Take $k$-points $t_1, \ldots , t_k \in S^1$ which are circularly ordered in $S^1$ 
(in particular, for the subindices of the $t_i$ we shall assume $i\in \ZZ/k\ZZ$ that is $k+1=1$).

\begin{figure}[ht]\begin{center}
\input{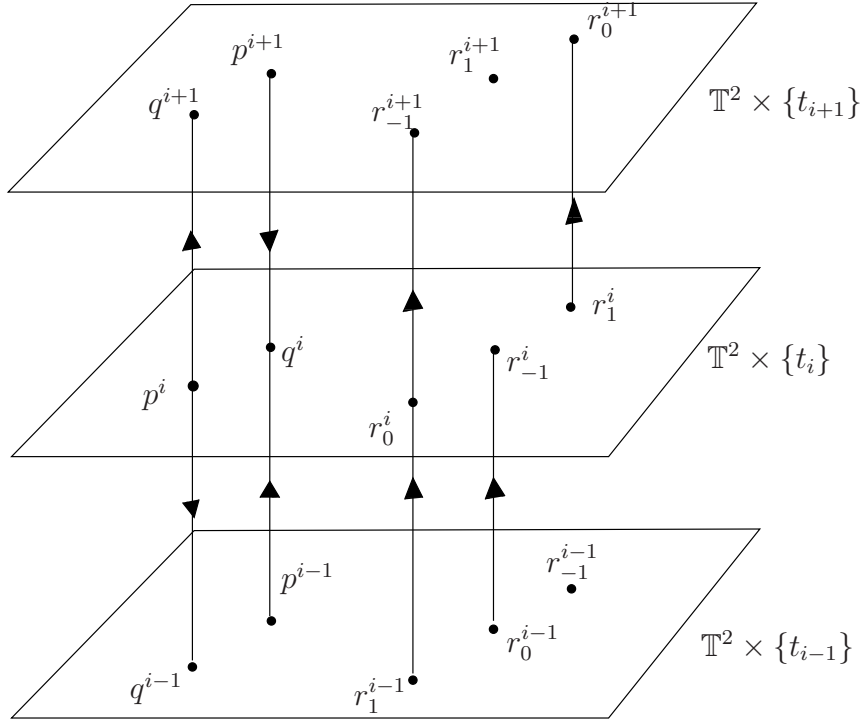}
\caption{\small{The picture of $f_0$ around $\TT^2 \times
\{t_i\}$}}\label{figure1}
\end{center}\end{figure}

We shall consider $f_0: \TT^3 \to \TT^3$ a $C^\infty$-diffeomorphism of the form $f_0(x,t)= (g_A(x), \varphi_x(t))$ 
satisfying the following properties: 

\begin{itemize}
\item[(P1)] $\varphi_x(t_i)=t_i$ for every $x\in \TT^2$ and every $t_i$; in other words, the torus $\TT^2\times \{t_i\}$ is invariant under $f_0$;

\item[(P2)]  $\frac 1 2 <
|(\varphi_x)'(t)| < \frac 3 2$ for every $(x,t) \in \TT^3$; 

\item[(P3)] in each torus $\TT^2 \times \{t_i\}$ we shall name $5$ fixed
points as $p^i=(\hat p^i,t_i)$, $q^i= (\hat q^i, t_i)$, $r^i_{-1}=(\hat r^i_{-1}, t_i)$, $r^i_{0}=(\hat r^i_{0}, t_i)$ and
$r^i_{1}=(\hat r^i_{1}, t_i)$. 

The points and the names are chosen so that
$$\hat p^i = \hat q^{i+1} = \hat p^{i+2}, \mbox{ and  }\hat r^{i}_{1}= \hat r^{i+1}_0 = \hat r^{i+2}_{-1}= \hat r^{i+3}_1,  \mbox{ for all }  i\in\ZZ/k\ZZ.$$

These fixed points  will satisfy that (see figure \ref{figure1}):
\begin{itemize}
\item[(a)] $|(\varphi_{\hat p^i})'(t_i)| >1$ and  $|(\varphi_{\hat
r^i_{1}})'(t_i)|>1$
,
\item[(b)] $|(\varphi_{\hat q^i})'(t_i)| <1$
and  $|(\varphi_{\hat r^i_{-1}})'(t_i)|<1$,

\end{itemize}
\item[(P4)] the following connections occur (see figure
\ref{figure1}):
\begin{itemize}
\item[(a)] every point in the interval 
$$\hat p^i \times (t_i,t_{i+1})= \hat
q^{i+1} \times (t_i,t_{i+1})$$ 
is contained in both  the unstable
set of $p^i$ and the stable set of $q^{i+1}$, 
\item[(b)] every
point in the interval  
$$\hat q^i \times (t_i,t_{i+1})= \hat p^{i+1} \times
(t_i,t_{i+1})$$
is contained in both  the stable set of $q^i$ and
the unstable set of $p^{i+1}$, 
\item[(c)] every point in the interval 
$$\hat
r^i_{0} \times (t_i, t_{i+1}) = \hat r^{i+1}_{-1} \times
(t_i,t_{i+1})$$
is contained in both  the unstable set of $r^i_0$
and the stable set of $r^{i+1}_{-1}$, 
\item[(d)] every point in the interval 
$$\hat r^i_{1} \times (t_i,t_{i+1}) = \hat r^{i+1}_{0} \times (t_i,
t_{i+1})$$
is contained in both the unstable  set of $r^i_1$ and
the stable set of $r^{i+1}_0$; 
\end{itemize}
\item[(P5)] there is $\nu>0$ so that, for every $t_i$ one has
that 
$$\log |(\varphi_x)'(t_i)| < - \nu$$
for every $x$ outside
small neighborhoods of the points $\hat p^i$, $\hat r^i_1$ and
$\hat r^i_0$ where $\log |(\varphi_x)'(t_i)| \leq
\frac{\nu}{2}$.
\end{itemize}

\begin{obs}\label{remark-odd} Property (P3) states that $\hat p^i = \hat q^{i+1} = \hat p^{i+2}$ as well as $\hat r^{i}_{1}= \hat r^{i+1}_0 = \hat r^{i+2}_{-1}= \hat r^{i+3}_1$.
It is for this to make sense in the complete order that we need that $k$ is divisible by $6$. 
 If $k$ were not divisible by $6$, simple modifications can be made to have a similar picture: For example, one can start with more fixed points and make more connections between the tori. 
 This will be evident in the proof.
\end{obs}

We have the following property which we will not use but shows why we say that these are many Kan's examples glued together.

\begin{prop}\label{Proposition-PropertiesF0}
The diffeomorphism $f_0$ is strongly partially hyperbolic with
splitting  $T\TT^3=E^s \oplus E^c \oplus E^u$ where $E^c$ is the
direction tangent to $\{x\} \times S^1$ and has exactly $k$ SRB
measures $\mu_1, \ldots, \mu_k$ such that $\mu_i$ is supported in
$\TT^2 \times \{t_i\}$. The union of their basins has full
Lebesgue measure and the basin of $\mu_i$ has positive measure in
every open subset of $\TT^2 \times [t_{i-1},t_{i+1}]$.
\end{prop}

\dem The fact that $f_0$ is strongly partially hyperbolic follows
directly from condition (P2) and the fact that $g_A$ is strongly
Anosov in the directions transverse to the circle fibers. For the claim about the SRB
measures, the proof is exactly the same as in Proposition 11.1 of
\cite{BDV} by using properties (P1)-(P5). See \cite{Kan,BV} or \cite{VY} for more details. \lqqd

To end this section, we note that $f_0$ preserves a center-stable
foliation $\cW^{cs}$ which consists of the leafs of the foliation
$W^s_A \times S^1$. As usual, we denote by $\cW^{cs}(x)$ the
leaf of the foliation through the point $x$.

We shall denote by $W^s_A (x,r)$ the interval of $W^s_A(x)$
centred at $x$ of length $r$.

\section{The first modification}\label{s.first}

We will consider a first modification $f_1: \TT^3 \to \TT^3$ obtained by making modifications close to each $\TT^2 \times \{t_i\}$ 
around the points $r^i_0$ of DA type (see for example \cite{BV}). 

We continue with the notation of the previous section.

Let $\eps>0$ be much smaller than the distance between any pair of tori $\TT^2 \times \{t_i\}$ and much smaller than the distance between any of the points $\hat p^i, \hat q^i, \hat r^i_0, \hat r^i_1, \hat r^i_{-1}$.

We will consider a smooth diffeomorphism $f_1: \TT^3 \to \TT^3$ which coincides with $f_0$ outside of $\bigcup_{i=1}^k B_{\eps}(r^i_0)$ and has the following properties:

\begin{itemize}
\item[(M1)] the foliation $\cW^{cs}$ is preserved by $f_1$;
\item[(M2)] the tori $\TT^2 \times \{t_i\}$ are preserved by $f_1$ and the dynamics along $\{r^i_0\} \times S^1$ remains unchanged for all $i$,
\item[(M3)] there exists a $Df_1$-invariant cone-field transverse to $\cW^{cs}$ and vectors in this cone-field are expanded by a factor of at least $3$ by $Df_1$. The cone-field is narrow enough so that:
\begin{itemize}
\item[(a)] there exists $L>0$ such that for every $x \in \TT^2$ and $\gamma$ a curve tangent to the cone-field of length $\geq L$ one has that $(W^s_{A}(x, 2\eps) \times S^1 )\cap \gamma \neq \emptyset$, 
\item[(b)] there is no curve tangent to the cone-field of length smaller than $L$ joining tori $\TT^2 \times \{t_i +\eps\}$ with $\TT^2 \times \{t_{i+1} - \eps\}$ for any $i$;

\end{itemize}
\item[(M4)] the maximal invariant set in $B_\eps(r^i_0) \cap (\TT^2 \times \{t_i\})$ is the point $r^i_0$ which is a hyperbolic source restricted to $\TT^2 \times \{t_i\}$ and two segments connecting $r^i_0$ to two other fixed points which are hyperbolic saddles restricted to $\TT^2 \times \{t_i\}$. These two hyperbolic saddles are not contained in $B_{\eps/2}(r^i_0)$. 

\end{itemize}

The DA-attractors contained in $\TT^2 \times \{t_i\}$ will be called $\Lambda_i$ and are characterised as the maximal invariant sets of $f_1$ in $(\TT^2 \times \{t_i\}) \setminus B_{\eps/2}(r^i_0)$, i.e. the complement of the basin of repulsion of $r_0^i$ in $\TT^2 \times \{t_i\}$.

The fact that such modifications can be made is quite classical (see \cite[Chapter 7]{BDV} or \cite[Section 6]{BV} for similar constructions).

One important feature  is the following:

\begin{prop}\label{Prop-MostlyContracting} For $f_1$, there exists a unique Gibbs-u-state $\mu_i$ supported on $\Lambda_i$ 
which is a SRB measure and has negative center Lyapunov exponents. 
\end{prop}

\dem The dynamics of $f_1$ in restriction to each torus $\TT^2\times \{t_i\}$ is Axiom A plus strong transversality, and admits a unique hyperbolic attractor. Therefore, the torus  $\TT^2\times \{t_i\}$
supports a unique Gibbs-u-state supported on $\Lambda_i$, which is indeed the unique SRB measure of the restriction.  It remains to see that this Gibbs-u-state is an SRB measure in the ambient manifold $\TT^3$. 
For that, it is enough to see that its transverse Lyapunov exponent is negative. 

Using property (P5) one can see that for every $x\in \Lambda_i$ one has that there is a positive Lebesgue measure subset of any 
unstable disk around $x$ such that the central Lyapunov exponent is negative (the argument is the same as in \cite{BV}, done for Ma\~n\'e's example, 
each unstable manifold has a full Lebesgue measure subset of points whose orbit spend most of the time out of the regions where the function $\varphi_x'(t_i)$ is $\geq 1$ so the average is negative). 
This concludes the proof, exactly as in \cite[Theorem B]{BV}. \lqqd

It is also rather simple to obtain that $\{q^1, \ldots, q^k \}$ is a skeleton. We will use the following result which will be important later too: 

\begin{lema}\label{l.stableqi} The stable manifold of $q^i$ contains $W^s_A(\hat q^i,3\eps) \times (t_{i-1},t_{i+1})$. 
\end{lema}
\dem  For $f_0$, the stable manifold of $q^i$ is exactly $W^s_A(\hat q^i) \times (t_{i-1},t_{i+1})$. Moreover, the set $W^s_A(\hat q^i,3\eps) \times (t_{i-1},t_{i+1})$ is forward invariant and disjoint from $\bigcup_{j=1}^k B_{\eps}(r^j_0)$, so, as $f_1$ coincides with $f_0$ outside $\bigcup_{j=1}^k B_\eps (r^j_0)$, one obtains the statement. \lqqd

We then obtain:

\begin{prop}\label{prop-skeletonf1} The diffeomorphism $f_1$ is mostly contracting and the set $\{q^1, \ldots, q^k\}$ is a skeleton.
\end{prop}

\dem To show that property (S1) is verified, it is enough to show that every unstable manifold intersects the stable manifold of one of the $q^i$. Indeed, if $x \in \TT^2 \times \{t_i\}$ it is immediate that $\cW^u(x)$ intersects the stable manifold of $q^i$ as $\Lambda_i$ is the unique attractor inside $\TT^2 \times \{t_i\}$ and it is hyperbolic. 

If a point $x$ belongs to $\TT^2 \times (t_i,t_{i+1})$ then it follows immediately from (M2) and Lemma \ref{l.stableqi} that $\cW^u(x)$ intersects both the stable manifold of $q^i$ and the stable manifold of $q^{i+1}$. This establishes (S1). 

Using Remark \ref{rem.mostlycontracting} and Proposition \ref{Prop-MostlyContracting} we get that $f_1$ is mostly contracting. 

It remains to show that  $\{q^1, \ldots , q^k\}$ verify property (S2), but this is immediate from the fact that the unstable manifold of $q_i$ is entirely contained in $\TT^2 \times \{t_i\}$. 
\lqqd

\begin{obs}
It is easy to show that the stable manifold of $q^i$ is dense in $\TT^2 \times (t_{i-1},t_{i+1})$. Indeed, given an open set $U \en \TT^2 \times (t_{i-1},t_{i+1})$, one can take an unstable arc $I \en U$  and after forward iteration, using property (M3) and the fact that the set $\TT^2 \times (t_{i-1},t_{i+1})$ is $f_1$-invariant one obtains that $f_1^n(I)$ intersects $W^s_A(\hat q^i,3\eps) \times (t_{i-1},t_{i+1})$ which is contained in $W^s(q^i)$ by Lemma \ref{l.stableqi}. Due to invariance, it holds that the closure of the stable manifold of $q^i$ is therefore exactly $\TT^2 \times [t_{i-1},t_{i+1}]$ and this is why we need to make a further perturbation to obtain the desired example. 
\end{obs} 


\section{Second modification, a small smooth perturbation}\label{s.second}

We will now make a very small smooth perturbation of $f_1$ in order to break the invariance of the tori $\TT^2 \times \{t_i\}$ so that having all basins intermingled becomes possible. 

The perturbation is made as follows. We consider $f: \TT^3 \to \TT^3$ which is $C^r$-close to $f_1$ and is obtained by composing $f_1$ with a diffeomorphisms $\rho^j$ supported on $B_{\eps/2}(r^j_0)$. 

The diffeomorphisms $\rho^j$ are the time $1$ maps of the flow of a vector field\footnote{We consider $\frac{\partial}{\partial t}$ to be the unit vector field tangent to the circles in $\TT^3 =\TT^2 \times S^1$.} $X= \delta \eta^j  \frac{\partial}{\partial t}$ where $\eta^j: M \to [0,1]$ is a smooth bump function which equals zero outside $B_{\eps/2}(r^j_0)$ and $1$ in $B_{\eps/4}(r^j_0)$. Notice that $\rho^j$ is arbitrarily $C^r$ close to the identity as $\delta \to 0$. 

\begin{figure}[ht]\begin{center}
\input{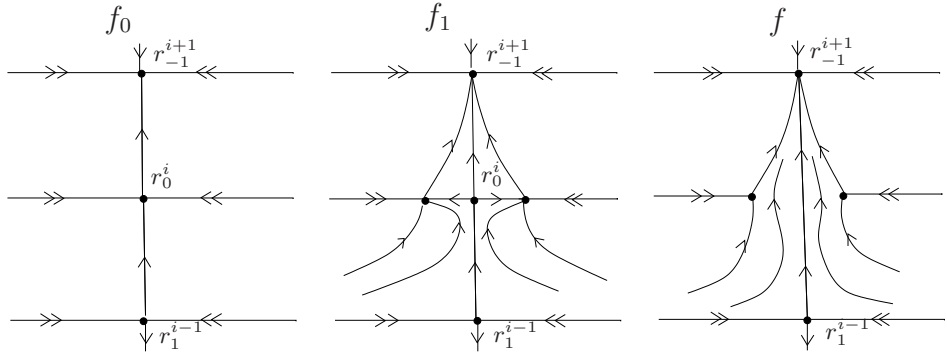}
\caption{\small{The dynamics in $\cW^{cs}(r^i_0)$. The modification in the $W^s_A \times S^1$ manifold in a neighborhood of $r^i_0$.}}\label{figure2}
\end{center}\end{figure}

Let us make some remarks about the properties of $f$ if $\delta$ is small enough: 

\begin{itemize}
\item[(R1)] conditions (M1) and (M3) of $f_1$ still hold for $f$ as the perturbation preserves the leafs of $\cW^{cs}$ and it is $C^1$-small;
\item[(R2)] $f$ is mostly contracting and $\{q^1, \ldots, q^k\}$ is a skeleton. This holds because the sets $\Lambda_i$ are disjoint from the support of the perturbation and therefore their $u$-Gibbs states remain with negative Lyapunov exponents, property (S1) is $C^1$-open and property (S2) still holds as the sets $\Lambda_i$, which coincide with the closure of the unstable manifold of $q^i$ remain unchanged. Indeed, Lemma \ref{l.stableqi} still holds and the proof of Proposition \ref{prop-skeletonf1} works exactly the same;
\item[(R3)] no point remains in $B_{\eps/2}(r^i_0)$. (See figure \ref{figure2}.)
\end{itemize}

The diffeomorphism $f$ is the diffeomorphism announced in Theorem \ref{teo-main}. By property (R2) it is enough to show that the stable manifolds of $q^i$ are dense in $\TT^3$ for all $i$ in order to be in the conditions of Theorem \ref{teo-basins}, this is done in Section \ref{SectionProof}. 

We remark that $f$ is partially hyperbolic with splitting $T\TT^3= E^{cs} \oplus E^u$ and \emph{dynamically coherent} (the foliation $\cW^{cs}$ is $f$-invariant and tangent to $E^{cs}$) but it is not strongly partially hyperbolic since in $\cW^{cs}(r^i_0)$ there are points in the stable manifold of $r^{i+1}_{-1}$ for
which vectors ``candidate'' to $E^s$ get mapped transverse to the strong stable of $r^{i+1}_{-1}$. Another reason for which $f$ cannot be strongly partially hyperbolic can be found in \cite{UV}. 


\section{Density of the stable manifold}\label{SectionProof}
To conclude the proof of the main theorem, we need to show that the stable manifold of each $q^i$ is dense which is the content of Proposition \ref{Prop-WsQiDense} below. 
For this, we need first a definition and some preliminary lemmas. 

If a hyperbolic saddle periodic point of $f$ has two contracting eigenvalues of different moduli, then one can  define both stable manifold $W^s(x)$ and  \emph{strong stable manifold} $W^{ss}(x)$ 
corresponding respectively to the eigenvalues of moduli smaller than one and of smallest modulus. The strong stable manifold of $x$ cuts  the stable manifold in two
components, called \emph{stable separatrices} of $x$.

If $x$ belongs to one of the $\Lambda_i\subset \TT^2\times \{t_i\}$, then its local stable and 
strong stable manifolds for $f$ coincide with the ones for $f_1$. Thus the local strong stable 
manifold of $x$ is contained in the torus $\TT^2\times \{t_i\}$ and the local stable manifold is a $2$-disc transverse to $\TT^2\times \{t_i\}$. Thus the transverse orientation of the tori 
$\TT^2\times \{t_i\}\subset \TT^2\times \SS^1$ (given by the orientation of the circle $\SS^1$) allows us to denote naturally by $W^s_+(x)$ and $W^s_-(x)$ the two stable separatrices of $x\in\Lambda_i$.



Notice that the boundary of the stable manifold  $W^s(q^i)$ in the center stable leaf  $\cW^{cs}(q^i)$ (in the intrinsic topology of $\cW^{cs}(q^i)$) 
contains the stable manifolds of $p^{i-1}$ and $p^{i+1}$ (which are one-dimensional). With our convention, the \emph{upper} $W^s_+(q^i)$  (resp. \emph{lower}  $W^s_-(q^i)$) separatrix of $q^i$  accumulates 
(in the topology of $\cW^{cs}(q^i)$) on $W^s(p^{i+1})$ (resp.  $W^s(p^{i-1})$).  

We now state some lemmas that will be needed.

\begin{lema}\label{lema-Wsqi} 
For all $i \in \{1, \ldots, k\}$ one has that\footnote{Here, one must understand $[t_i,t_{i-1})$ as $(t_{i-1},t_i]$. }
\begin{equation}\label{eq:halfmfd} 
 W^s_{\pm} (q^i) = \bigcup_{n\geq 0} f^{-n}(W^s_A(\hat q^i, 3\eps) \times [t_i, t_{i\pm1}))  \ . 
 \end{equation}
\end{lema}

\dem One notices that the set $\hat W= W^s_A(\hat q^i, 3\eps) \times (t_{i-1},t_{i+1})$ is contained in $W^s(q^i)$ for $f_0$ and inside there, $W^s_A(\hat q^i, 3\eps) \times \{t_i\}$ 
corresponds to the intersection with the strong stable manifold for $f_0$. Since the forward iterates of the points in $\hat W$ under $f_0$ remain in $\hat W$ and $\hat W$ 
is disjoint from the region perturbed to obtain $f_1$ and $f$, one deduces that $\hat W$ is contained in the stable manifold of $q^i$ for $f$ and that $W^s_A(\hat q^i, 3\eps) \times \{t_i\}$ 
contains a fundamental domain of the strong stable manifold. Now it is clear that iterating backwards one obtains the complete invariant manifolds and since orientation is preserved 
one concludes 
that the half spaces also have the desired property.
\lqqd

\begin{figure}[ht]\begin{center}
\input{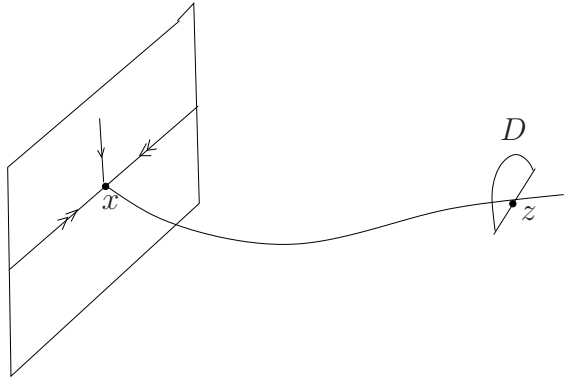}
\caption{\small{The hypothesis of Lemma \ref{lema-lambda}.}}\label{figure3}
\end{center}\end{figure}

The following is a restatement of the classical $\lambda$-Lemma for manifolds with boundary in our context (see figure \ref{figure3}). 

\begin{lema}\label{lema-lambda} 
Let $x$ be a hyperbolic fixed point of $f$ with two eigenvalues of different modulus and both smaller than one. Let $D$ be a half disk\footnote{i.e. the image of an embedding of $\{(t,s) \in \RR^2 \ : \ t^2 +s^2 <1  \ , \ s\geq 0\}$.} so that 
\begin{itemize}
\item there exists $z \in \partial D \cap W^u(x)$ such that $T_z D \oplus W^u(x) = T_z \TT^3$ ,
\item  the limit of $Df^{-n}(T_z\partial D)$ converges as $n\to \infty$ to the eigenspace of the smallest modulus of $D_xf$. 
\end{itemize}
Then one has that $f^{-n}(D)$ converges on compact sets to either $W^s_+(x)$ or $W^s_-(x)$ depending on the orientation of $D$ with respect to $z \in \partial D$.  
\end{lema}

\dem This is a very standard result. We refer the reader to e.g \cite[Page 4860]{BCGP} for a very similar statement with proof. 
\lqqd

\begin{obs}\label{rem-clasiclamda}
The classical $\lambda$-Lemma implies that if $z$ is in the interior of $D$ then $f^{-n}(D)$ converges in compact parts to the whole $W^s(x)$. 
\end{obs}

\begin{lema}\label{lem-induction}
The closure of  $W^s_-(q^i)$ in $\TT^3$ contains $W^s_-(q^{i-1})$ and $W^s_+(q^{i-2})$. 
\end{lema}

\dem Since $W^{ss}(q^i) \cap \cW^u(r^i_{-1}) \neq \emptyset$ one can apply Lemma \ref{lema-lambda} to obtain that $W^{s}_-(r^i_{-1}) \subset \overline{ W^s_-(q^i)}$, notice that $W^{ss}(q^i)$ divides a small lower disk of the stable manifold $W^s(q^i)$ and the backward iterates of its tangent space converges to the eigenspace of smallest eigenvalue of $D_{r^i_{-1}}f$ since the set  $\Lambda_i$ is strongly partially hyperbolic. As the set $W^s_-(q^i)$ is invariant, we deduce the statement.  

Notice that inside $\cW^{cs}(r^i_1)$ the stable manifold of $r^i_{-1}$ contains an arc transverse to the unstable manifold of $r^{i-2}_1$, and therefore, by backward iteration it contains a band above the (strong) stable manifold of $r^{i-2}_1$ by the standard $\lambda$-lemma. 

Now, notice that the strong unstable manifold of $q^{i-2}$ intersects the stable manifold of $r^{i-2}_1$ and one can therefore apply Lemma \ref{lema-lambda} to deduce that $W^s_+(q^{i-2}) \en \overline{W^s_-(r^i_{-1})} \en \overline{W^s_-(q^i)}$.

\begin{figure}[ht]\begin{center}
\input{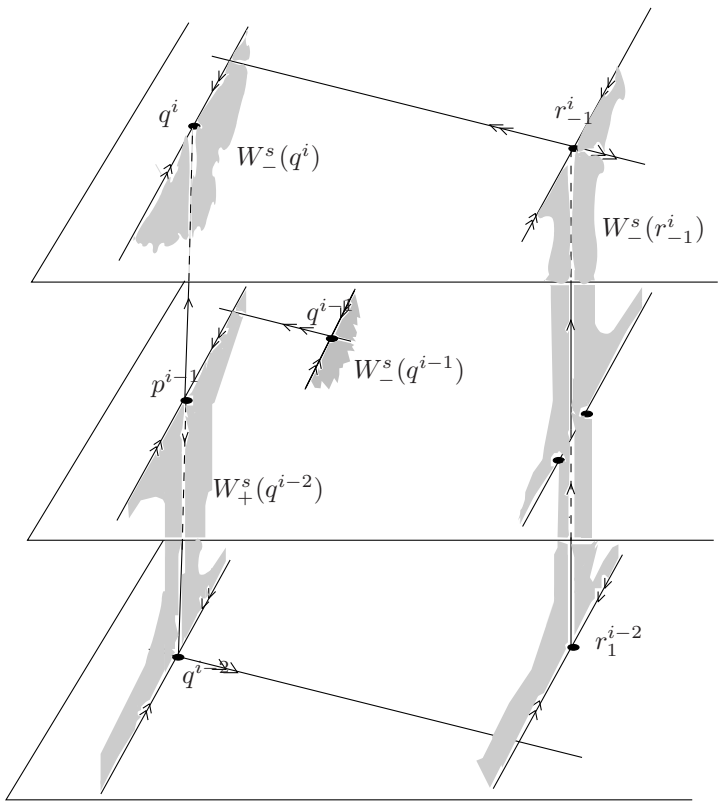}
\caption{\small{The argument in Lemma \ref{lem-induction}. $W^s_{-}(q^{i})$ accumulates on $W^s_{-}(r^i_{-1})$ which in turn accumulates on $W^s_+(q^{i-2})$. As $W^s_+(q^{i-2})$ accumulates on $W^s(p^{i-1})$ from below, it accumulates on $W^s_-(q^{i-1})$. }}\label{figure4}
\end{center}\end{figure}

To conclude, we apply Lemma \ref{lema-Wsqi} to $q^{i-2}$ and using an intersection of $\cW^u(q^{i-1})$ with $W^s_A(\hat q^{i-2}, 3\eps) \times [t_{i-2}, t_{i-1}))$ we obtain that $W^s_-(q^{i-1}) \en \overline{ W^s_+ (q^{i-2})}  \en \overline{ W^s_-(q^i)}$.  This concludes the proof of the lemma  (see figure \ref{figure4}). \lqqd

\begin{prop}\label{Prop-WsQiDense}
The stable  manifold of $q^i$ is dense in $\TT^3$ for every $i=1,\ldots, k$.
\end{prop}

\dem  Applying Lemma \ref{lem-induction} for all $i$ and invariance one obtains that the closures of all stable manifolds of $q^i$ all coincide and therefore are dense (e.g. as $\{q^1, \ldots, q^k \}$ is a skeleton). 
\lqqd

This implies that all basins are intermingled by Theorem \ref{teo-basins} and completes the proof of Theorem \ref{teo-main}.

\section{Countably many SRB measures with intermingled basins}\label{sec-countable}

In this section we will explain how some modifications to the construction above allow us to construct an example with countably many intermingled basins. 

For this to be possible, one will start with countably many tori accumulating in a given one from one side. We denote, as before $f_0: \TT^3 \to \TT^3$ to such a diffeomorphism and we assume the tori are $T_n = \TT^2 \times \{1/n\}$ and $T_\infty=\TT^2 \times \{0\}$ the limit torus where we are considering $S^1 = [0,1]/_{0 \sim 1}$ and $n\geq 2$. 

For the initial dynamics, every tori $T_n$ will support an SRB measure $\mu_n$ (mostly contracting) except the limit one. Notice that one is forced to make the Lyapunov exponents in each such measure to go to $0$ as the limit measure $\mu_\infty$ in $T_\infty$ cannot have negative nor positive center Lyapunov exponents. As we will not be able to "cross" the torus $T_\infty$ we will be forced to make connections going "up" and "down". This means that conditions (P1)-(P5) will have to be adapted accordingly to have at least 7 fixed points and construct in each torus a fixed point which in the center:

\begin{itemize}
\item is attracting  ($q^i$)
\item  one which is repelling ($p^i$)
\item one saddle node going up (as the $r^i_0$) but also one going down (call them $s^i_0$) which have corresponding $s^{i+1}_{1}$ and $s^{i-1}_{-1}$. 
\end{itemize}

Now, the construction of the DA deformation to obtain $f_1$ has to be done in two steps, first, one has to slow down the stable exponents of all the fixed points in $T_\infty$ which will be limits of $r^i_0$ and $s^i_0$ so that the modifications can be done without loosing the regularity. 

Now, as the DA-diffeomorphisms are $C^\infty$ close to diffeomorphisms of tori having points with zero stable exponents in the fixed points, one can construct $f_1$ by opening in each $r^i_0$ and $s^i_0$ but of course the amount one opens goes to zero as $i \to \infty$. By 'opening in each $r^i_0$ and $s^i_0$' we mean that we make the same construction as in Section \ref{s.first}, see figure \ref{figure2}.

The same statements we have obtained in Section \ref{s.first} adapted accordingly hold for this $f_1$, there are two points to be careful: 

\begin{itemize}
\item The limit torus $T_\infty$ admits a Gibbs-u-state which is not an SRB measure. In particular, the sequence $\{q^2, \ldots, q^n, \ldots \}$ is not a skeleton. However, their stable manifolds intersect the unstable manifold of every point which does not intersect $T_\infty$. 

\item Not being a skeleton, the property is no longer robust, but it is if one makes a $C^1$ perturbation which is the identity in $T_\infty$. 
\end{itemize}

\begin{figure}[ht]\begin{center}
\input{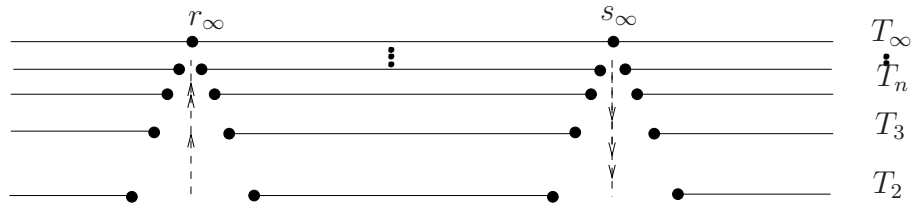}
\caption{\small{Schematic figure of the accumulating tori and openings getting smaller.}}\label{figure5}
\end{center}\end{figure}

Finally, one makes a small $C^\infty$-perturbation $f$ of $f_1$ which pushes up in the middle of the $r^i_0$ and down in the $s^i_0$. This is the desired example and properties of Section \ref{s.second} hold accordingly, see figure \ref{figure5}.

To see that all basins are intermingled, again, one has to use an adaptation of Theorem \ref{teo-basins} as the set of points $\{q^n\}_n$ is not really a skeleton. But the proofs adapt directly (as one shows that the basins are pairwise intermingled). Showing the density of the stable manifolds of the points $q^i$ requires a bit more work as one has points moving up and down, but a more careful use of Lemmas \ref{lema-Wsqi} and \ref{lema-lambda} give that: 

\begin{lema} 
The closure of $W^s_-(q^i)$ of $q^i$ contains $W^s_+(q^{i-2})$ and $ W^s_{-}(q^{i-1})$; and the closure of $W^s_+(q^i)$ contains $W^s_-(q^{i+2})$ and $W^s_+(q^{i+1})$. 
\end{lema}

Putting all this together, one obtains an example with countably many intermingled basins. 
\section{Volume hyperbolicity}\label{sec-volumehyperbolic} 

In this section we will show how one can ensure that the splitting $T\TT^3 = E^{cs} \oplus E^u$ of $f$ is volume hyperbolic by making some of the choices in the construction of $f_1$ and $f$ more precise. This will help proving the topological mixing of a specific example in the next section. 

One must add a new property to the modification made for $f_1$ which involves some constants that we will choose appropriately when constructing $f$. 

\begin{itemize}
\item[(M5)] There are constants $\xi,\zeta>0$ such that outside the $\zeta$-neighbourhood of the points $r^i_0$ the volume along $E^{cs}$ is contracted by $Df_1$ by a factor $1/2$ and the maximum dilatation of volume of $Df_1$ along $E^{cs}$ in all of $\TT^3$ is of $1+\xi$. 
\end{itemize}

Notice that for $f_0$ one has volume contraction along $E^{cs}$ of a factor at least $1/3$ as the contraction along $E^s$ is smaller than $1/5$ and the expansion along $E^c$ smaller than $3/2$ and these are orthogonal. As in a neighbourhood of $r^i_0$ for $f_1$ there is some half disc which is mapped outside itself, one needs some amount of expansion, which we are controlling with condition (M5) to be not so large and concentrated in a small neighbourhood of $r^i_0$. The fact that the area expansion can be bounded follows from the fact that for $f_0$ the derivative of $r^i_0$ along $E^c$ equals one as it is a saddle node point. 

An important remark is that, if one fixes $\zeta\ll \eps$, the size of a perturbation of $f_1$ which will verify conditions (R1)-(R3) can be chosen independent of the constant $\xi$. Using this remark, one can choose $\delta$ (for the perturbation of $f_1$ giving rise to $f$) in order that the number of iterates that a point remains in $B_\zeta(r^i_0)$ is uniformly bounded by a constant $N_0$ which is independent of $\xi$. Now, let $N_1>0$ be such that $f^k(B_\zeta(r^i_0)) \cap B_\zeta(r^j_0) = \emptyset$ if $i\neq j$ and $1\leq k \leq N_1$, notice that as $\zeta \ll \eps$ such an $N_1\geq 2$ exists and it is independent of $\xi$ too. The perturbation made to create $f$ from $f_1$ can be chosen to be a translation in each $B_\zeta(r^i_0)$ for every $i$ and therefore does not change the rate of contraction of area along $E^{cs}$. 

It follows that if one chooses $\xi$ small enough so that $(1+\xi)^{N_0} (\frac {1}{2})^{N_1} < 1$ one obtains that $Df|_{E^{cs}}$ contracts the area uniformly and therefore $f$ is volume hyperbolic. 

\begin{obs}\label{obs-volumecountable}
Notice that volume hyperbolicity is not possible in our construction of an example with countably many intermingled basins as the limit torus should have some periodic points with jacobian equal to one along the $cs$-direction. However, one can obtain volume hyperbolicity out of arbitrarily small neighbourhoods of these points. 
\end{obs}

\section{Topologically mixing examples}\label{sec-recurrence}

In this section we will indicate some variants on the construction that allow one to obtain examples which are \emph{topologically mixing}. We will use some classical results on robust heterodimensional cycles and blenders for which we refer the reader to \cite[Chapters 6 and 7]{BDV} for a broad introduction. To obtain this property, we need to introduce some additional properties to the ones presented in Section \ref{s.Kan} to $f_0$ in order to guarantee those properties. 

These are: 

\begin{itemize}
\item[(P6)] The diffeomorphism $g_A$ has two additional fixed points, $\hat a$ and $\hat b$. 

\item[(P7)] In a small neighbourhood of $\hat a$ the map $\varphi_x: S^1 \to S^1$ restricted to $[t_1,t_2]$ has exactly three hyperbolic fixed points, $t_1$, $t_2$ which are attracting and  a point $ t_a \in (t_1,t_2)$ which is repelling . We call the fixed point $a=(\hat a, t_a)$ of $f_0$. 

\item[(P8)] In a small neighbourhood of $\hat b$ the map $\varphi_x: S^1 \to S^1$ restricted to $[t_1,t_2]$ has exactly three hyperbolic fixed points, $t_1$, $t_2$ which are repelling\footnote{This means that one has to adapt (P5) to take this into account.} and  a point $ t_b \in (t_1,t_2)$ which is attracting . We call the fixed point $b=(\hat b, t_b)$ of $f_0$. 

\item[(P9)] The points $A$ and $B$ are connected by a robust cycle and the closure of their stable and unstable manifolds coincide robustly. 
\end{itemize}

There are many ways to construct examples with this properties. We refer the reader to \cite{BD,BDK} for similar constructions. The key feature is the use of \emph{blenders} which we will not define here. We will also choose the construction of $f$ so that it is volume hyperbolic as explained in Section \ref{sec-volumehyperbolic}. 

With this additional properties, one can show:

\begin{prop} The diffeomorphism $f$ is topologically mixing.
\end{prop}

\dem As a direct consequence of the $\lambda$-Lemma we know that if there is a hyperbolic fixed point $P$ for which both the stable and unstable manifold are dense, then $f$ is topologically mixing. 

We need to show that the stable and unstable manifold of, say, $a$ are dense. Using property (P9) (which persists after the modifications) it is enough to show that the stable manifold of $b$ as well as the unstable manifold of $a$ are dense. 

The fact that the stable manifold of $b$ is dense follows exactly as for the points $q^i$ (see Proposition \ref{Prop-WsQiDense}). 

The density of the unstable manifold of $a$ is more involved and uses the volume hyperbolicity of $f$. 

First, we show that the union of the unstable manifolds of the points $p^i$ is dense in $\TT^3$. Notice that this union contains the unstable manifold of the circle $\{\hat p^i\} \times S^1$ which is normally hyperbolic, denote as $W^u_{loc}=W^u_A(\hat p^i, 3\eps) \times S^1$ its local unstable manifold (as there are no perturbations in this region, see Lemma \ref{lema-Wsqi} for an identical argument). As $\cW^{cs}$ remains unchanged and is minimal, it follows that for each point $x\in \TT^3$ the center stable leaf $\cW^{cs}(x)$ is intersected by $W^u_{loc}$ in two circles leaving a bounded cylinder containing $x$ and with uniformly bounded area. Using volume hyperbolicity, this implies that the forward iterates of $W^u_{loc}$ become dense proving our claim. 

Now, we show that all the $p^i$ are homoclinically related, therefore, their unstable manifolds are all dense. Indeed, each $p^i$ is homoclinically related with $r^i_1$ and the stable manifold of $p^{i+1}$ intersects the unstable manifold of $r^{i}_1$ transversally, therefore the stable manifold of $p^{i+1}$ intersects the unstable manifold of $p^i$ transversally. Inductively, one obtains that all $p^i$ are homoclinically related.  As a byproduct of this argument, using remark \ref{rem-clasiclamda} one obtains that the positive unstable separatrix of $r^i_1$ is dense in $\TT^3$ for every $i$ as the stable manifold of $p^i$ intersects its interior transversally and its unstable manifold is dense. Notice that the \emph{positive unstable separatrix} is well defined in analogy to the ones defined in Section \ref{SectionProof}.

To conclude, it is therefore enough to show that the stable manifold of some $r^i_1$ intersects the unstable manifold of $a$ with a configuration as in Lemma \ref{lema-lambda}. Recall that that the unstable manifold of $a$ consists of the backward iterates of $W^u_A(\hat a, 3\eps) \times (t_1,t_2)$. Therefore, one gets that the strong stable manifold of $r^1_0$ intersects the boundary of $W^u(a)$ (relative to the unstable manifold of the normally hyperbolic circle $\{\hat a\} \times S^1$) and one can apply Lemma \ref{lema-lambda} to get density of the unstable manifold of $a$ as it will become dense in the positive half unstable manifold of $r^1_0$. This concludes the proof of the proposition. 
\lqqd

\section{Final comments}

We point out that to make the examples with countably many intermingled basins topologically mixing one should need to adapt the arguments in the previous section as in principle, volume hyperbolicity is not available (c.f. Remark \ref{obs-volumecountable}). Using the fact that this volume hyperbolicity can be achieved in very large regions of $\mathbb{T}^3$ it is certainly likely that the arguments can be carried out, still, we pose this as a question:

\begin{quest}
Can the examples with countably many intermingled basins be made topologically mixing?
\end{quest}

Finally, we point out that it is also likely that the construction made in $\mathbb{T}^3$ using tori as the initial attractors can be extended to other hyperbolic attractors in surfaces such as the Plykin attractor. Doing this one may create for instance, a solid torus with Plykin attractors in the disk and a skew product dynamics similar to the one made here to get many intermingled basins in a solid torus. As the solid torus can be embedded in any 3-manifold this would give examples in any isotopy class of any 3-dimensional manifold. 





\begin{thebibliography}{2}

\bibitem[BCGP]{BCGP} C. Bonatti, S.Crovisier, N. Gourmelon, R. Potrie, Tame dynamics and robust transitivity: chain recurrence classes versus homoclinic classes, \emph{Transactions of the AMS} {\bf 366} (9) (2014) 4849--4871.

\bibitem[BD]{BD} C. Bonatti, L. D\'iaz, Persistent nonhyperbolic transitive diffeomorphisms, \emph{Ann. Math.} {\bf 143} (1995) 367?396.

\bibitem[BDK]{BDK} C. Bonatti, L. D\'iaz, S. Kiriki, Stabilization of heterodimensional cycles, \emph{Nonlinearity} {\bf 25} (2012) 931--960.

\bibitem[BDV]{BDV} C. Bonatti, L. D\'iaz, M. Viana, \emph{Dynamics Beyond Uniform Hyperbolicity. A global geometric and probabilistic perspective}, Encyclopaedia of Mathematical Sciences {\bf 102}. Mathematical Physics III. Springer-Verlag (2005).

\bibitem[BV]{BV} C. Bonatti and M. Viana, SRB measures for partially hyperbolic diffeomorphisms whose central direction is mostly contracting. \emph{Israel J. of Math} {\bf 115} (2000), 157--193.

\bibitem[DVY]{VY} D. Dolgopyat, M. Viana, J. Yang, Geometric and measure-theoretical structures of maps with mostly contracting center,  to appear in  \emph{Comm. Math. Physics}  arXiv:1410.6308.

\bibitem[HHTU]{HHTU} F. Rodriguez Hertz, J. Rodriguez Hertz, A. Tahzibi, R. Ures, Uniqueness of SRB measures for transitive surface diffeomorphisms, \emph{Comm. Math. Physics} {\bf 306} (2011) 35--49.

\bibitem[K]{Kan} I. Kan, Open sets of diffeomorphisms having having two attractors, each with an everywhere dense basin, \emph{Bull. Amer. Math. Soc.} {\bf 31} (1994) 68-–74.

\bibitem[MW]{MW} I. Melbourne, A. Windsor,  A $C^\infty$ diffeomorphism with infinitely many intermingled basins, \emph{Ergodic theory and dynamical systems} {\bf 25} (2005) 1951--1959.

\bibitem[M]{Milnor} J. Milnor, On the concept of attractor, \emph{Comm. Math. Phyisics} {\bf 99} 2 (1985), 177-195.


\bibitem[UV]{UV} R. Ures, C. Vasquez, On the robustness of intermingled basins, \emph{Preprint}  arXiv:1503.07155.



\end{thebibliography}
\end{document}